\documentclass[a4paper, 10pt,oneside]{article}
\usepackage{amsfonts,amsmath, amssymb,amsthm}
\usepackage[T1]{fontenc}
\usepackage[latin9]{inputenc}
\usepackage{setspace}
\usepackage{amssymb}

\newtheorem{thm}{Theorem}[section]

\newtheorem{defn}[thm]{Definition}
\newtheorem{example}[thm]{Example}
\newcommand{\Rset}{\mathbb {R}}
\newcommand{\Cset}{\mathbb {C}}
\newcommand{\Zset}{\mathbb {Z}}

\newcommand{\RR}{\mathcal {R}}
\newcommand{\A}{\mathcal {A}}
\newcommand{\B}{\mathcal {B}}

\newcommand{\DD}{\mathcal {D}}

\newcommand{\F}{\mathcal {F}}

\newcommand{\dd}{\mathrm {d}}
\newcommand{\e}{\mathrm {e}}
\newcommand{\bfb}{\mathbf {b}}
\newcommand{\ii}{\mathrm {i}}
\newcommand{\rank}{\mathrm {rank}\,}

\newcommand{\Ker}{\mathrm {Ker}\,}

\renewcommand{\Re}{\mathrm{Re}}

\begin{document}
\title{Stability, stabilizability and exact controllability\\of a class of linear
neutral type systems}
\author{R. Rabah$^*$  \hskip 5ex G. M. Sklyar$^\dag$\\
{\normalsize $^*$IRCCyN/\'Ecole des Mines}\\
{\normalsize 4 rue Alfred Kastler, 44307 Nantes, France}\\
{\normalsize $^\dag$Institute of mathematics, University of Szczecin}\\
{\normalsize Wielkopolska 15,70451 Szczecin, Poland}\\
\tt{\small rabah@emn.fr} \hskip 2ex \tt{\small sklar@univ.szczecin.pl}
}
\date{}
\maketitle
\textbf{Keywords:} Neutral type systems, exact controllability, stabilizability, stability

\begin{abstract}
Linear systems of neutral type are considered using the infinite dimensional approach. The main problems are asymptotic, non-exponential stability, exact controllability and regular asymptotic stabilizability. The main tools are the moment problem approach, the Riesz basis of invariant subspaces and the Riesz basis of family of exponentials. 
\end{abstract}
\section{Introduction}
Many applied problems from physics, mechanics, biology,
and other fields can be described by delay differential equations.  A large class of such systems are systems of neutral type. In this paper, we consider  a general class of
 neutral systems with distributed delays given by the equation
\begin{equation}\label{1}
\left \{
 \begin{array}{l}
\displaystyle{\frac{\dd}{\dd t}}[z(t)-K z_t]=L z_t+Bu(t),\quad  t\ge 0,\cr
z_0=\varphi,
\end{array}
\right.
\end{equation}
where $z_t:[-1,0]\to \Cset^n$ is the history of $z$ defined by
$z_t(s)=z(t+s)$. The difference and delay operators $K$ and $L$,
respectively, are defined by
\begin{eqnarray*}
K f=A_{-1}f(-1)\quad\mathrm{and}\quad
L f=\int^0_{-1}A_2(\theta)\frac{\dd}{\dd \theta}f(\theta)\,\dd\theta+\int^0_{-1}A_3(\theta)f(\theta)\,\dd\theta
\end{eqnarray*}
for $f\in H^1([-1,0],\Cset^n)$,
where $A_{-1}$ is a constant $n\times n$ matrix,
 $A_2, A_3$ are $n\times n$ matrices whose
 elements belong to $L_2(-1,0)$, and $B$ is a constant $n\times r$ matrix.
 
A more general case may be when
$$
K f=\int^0_{-1}\,\dd \mu(\theta)f(\theta),\quad f\in C([-1,0],\Cset^n),
$$
where $\mu$ is a matrix-valued function of bounded variation and continuous at zero. But we limit ourself to
above mentioned case when $K f=A_{-1}f(-1)$ because that the property of system are mainly characterized by the structure of the matrix $A_{-1}$. For more general forms of the operator $K$ it may expected that the situation is analogous, but it is not clear how the properties studied here may be connected. Anyway, it is a domain for further investigation. Distributed delay may arise in the natural modeling or after some feedback.
Our purpose is to investigate the problems of asymptotic stability, of stabilizability by linear feedback and of exact controllability.

For that problems, the neutral type systems are less studied that the retarded systems when $K=0$. The difficulties are related to the following particular properties of neutral type systems: there may exist an infinite number of eigenvalues in the right half plane, in particular near the imaginary axis; the choice of the phase-space is crucial, in contrast to the case of retarded functional
differential equations where solutions are more smooth than the initial data; some feedback with may change the structure of the system, etc.
\section{The operator model}
In \cite{hale} and several other works, the framework is based on the
description of neutral type systems in the space of continuous
functions $C([-1,0]; {\Cset}^n)$.  However, for several control problem, the Hilbert space structure is more convenient in the study of our class of systems. In Hilbert spaces one can use the fundamental tool of Riesz basis (or orthonormal basis modulo a bounded isomorphism). We consider the operator model of neutral
type  systems introduced by  Burns and al. in product spaces.
This approach  was also used in \cite{lunel-yak} for the construction of a spectral model. 
In \cite{yama} the authors consider the particular case of discrete delay,
which served as a  model in \cite{rs-aml,rab_skl2} to characterize the stabilizability
of a class of systems of neutral type.
 
The state space is $M_2(-1,0; \Cset^n)=\Cset^n \times L_2(-1,0; \Cset^n)$, briefly $M_2$, and
permits (\ref{1}) to be  rewritten as
\begin{equation}\label{2s}
 \frac{{\mathrm d}}{{\mathrm d}t}x (t)= {\mathcal A} x(t)  + \B u(t), \qquad x(t)= \begin{pmatrix}{y(t)}\cr{z_t(\cdot)}\end{pmatrix},
\end{equation}
where the operators $\A$ and $\B$ are defined by
\begin{equation} \label{AB}
{\mathcal A}
 \begin{pmatrix}{y(t)}\cr { z_t(\cdot)}\end{pmatrix}=  \begin{pmatrix}{\int^0_{-1}A_2(\theta)\dot z_t(\theta){\mathrm
d}\theta  + \int^0_{-1}A_3(\theta) z_t(\theta){\mathrm
d}\theta }\cr { {\mathrm d}z_t(\theta)/ {\mathrm d}\theta}\end{pmatrix}, 
\qquad \B u=\begin{pmatrix} Bu\cr 0\end{pmatrix}
\end{equation}
The domain of ${\mathcal A}$ is given by 
$${\mathcal
D}({\mathcal A})= \{ (y,z(\cdot)) : z\in H^1(-1,0;\Cset^n),
y=z(0)-A_{-1}z(-1)\} \subset M_2
$$ 
and  the operator
${\mathcal A}$ is the infinitesimal generator of a
$C_0$-semigroup $\e^{{\mathcal A}t}$. The relation between the solutions 
of the delay system (\ref{1})
and the system (\ref{2s}) is $z_t(\theta)=z(t+\theta)$.

In the particular case when $A_2(\theta)=A_3(\theta)= 0,$
we use the notation $\widetilde{\mathcal A}$ for ${\mathcal
A}.$ The properties of $\widetilde{\mathcal A}$
can be expressed mainly in terms of the properties of matrix
$A_{-1}$ only. Some important properties of
${\mathcal A}$ are close to those of $\bar{\mathcal A}$.
\section{Spectral analysis}
Let us denote by $\mu_1,...,\mu_\ell$, $\mu_i\not=\mu_j$ if
$i\not=j$, the eigenvalues of $A_{-1}$ and the dimensions of their
rootspaces (generalized eigenspaces)
 by $p_1,...,p_\ell,$ $\sum_{k=1}^\ell p_k=n.$ Consider the points
$\lambda^{(k)}_m\equiv \ln |\mu_m| + i(\arg \mu_m + 2\pi
k), m=1,..,\ell; k\in Z $ and the circles $L^{(k)}_m$ of
fixed radius $r\le r_0\equiv \frac{1}{3}\min \{
|\lambda^{(k)}_m -\lambda^{(j)}_i |, (m,k)\neq (i,j)\}$
centered at $\lambda^{(k)}_m.$

\begin{thm}\label{spec}
 The spectrum of ${\cal A}$ consists of
the eigenvalues only which are the roots of the equation $\det
\Delta (\lambda)=0,$ where
\begin{equation}\label{Del}
\Delta_{\cal A} (\lambda)=\Delta (\lambda)\equiv -\lambda I +
\lambda e^{-\lambda}A_{-1} + \lambda   \int^0_{-1}e^{\lambda s}
A_2(s)\dd s  + \int^0_{-1}e^{\lambda s} A_3(s)\dd s.
\end{equation}
The corresponding eigenvectors of ${\cal A}$ are $\varphi= \left( 
C-e^{-\lambda}A_{-1}C, e^{\lambda\theta} C \right),$ with
$C\in \Ker  \Delta (\lambda).$

There exists $N_1$ such that for any  $|k|\ge
N_1,$ the total multiplicity of the roots of the equation
$\det \Delta (\lambda)=0,$ contained in the circle
$L^{(k)}_m,$ equals $p_m.$
\end{thm}

The  description of the location of the spectrum of ${\cal A}$ we use
Rouch\'e theorem.

\subsection{Basis of invariant subspaces}
The most desired situation for concrete systems is to have
a Riesz basis formed by eigenvectors of $A$ or, at least,
by generalized eigenvectors.
In more general situations, one studies the existence of
basises formed by subspaces. We remind that {\it a sequence
of nonzero subspaces $\{ V_k\}^\infty_i$ of the space $V$
is called basis (of subspaces) of the space $V$, if any
vector $x\in V$ can be uniquely presented as
$x=\sum^\infty_{k=1} x_k,$ where $x_k\in V_k, \, k=1,2,..$}
We say that the basis $\{ V_k\}^\infty_i$ is orthogonal if
$V_i$ is orthogonal to $V_j$ when $i\neq j.$  
A basis $\{ V_k\}$ of subspaces is called  a Riesz basis if there are an orthogonal basis of
subspaces $\{ W_k\}$  and a linear bounded invertible operator
$R,$ such that $RV_k = W_k.$

The best "candidates" to form the basis of subspaces are
generalized eigen\-spaces of the generator of a semigroup,
but there are simple examples (see Example~\ref{ex1} below) showing
that generalized eigenspaces do not form such a basis in
the general case.

One of the main ideas of our approach is to construct a Riesz
basis of finite-dimensional subspaces which are invariant for the
generator of the semigroup  (see (\ref{2s})).  

In \cite{rsr-CR,rsr-stab} we obtained the following general result.

\begin{thm}\label{th-basis}
 There exists a sequence of invariant for
${\cal A}$  (see (\ref{2s})) finite-dimensional subspaces which
constitute a Riesz basis in $M_2.$

More precisely, these subspaces are $\{ V^{(k)}_m, |k|\ge
N, m=1,..,\ell \}$ and a $2(N+1)n$-dimen\-sional subspace
spanned by all eigen- and rootvectors,
 corresponding to all eigenvalues of ${\cal A},$
 which are outside of all
circles $L^{(k)}_m,$ ${|k|\ge N, m=1,..,\ell}.$

Here $V^{(k)}_m \equiv P^{(k)}_m M_2$,
where 
$$P^{(k)}_m M_2=
\frac{1}{2\pi \ii}\int_{L^{(k)}_m} R({\cal A},\lambda) d\lambda
$$ 
are
spectral projectors; $L^{(k)}_m$ are circles defined before.
\end{thm}

We emphasize that the operator $ {\cal A}$ may not possess in
a Riesz basis of generalized eigenspaces. We
illustrate this on the following

\begin{example} \label{ex1}
{\em Consider the particular case of the
system (\ref{1}):
\begin{equation}\label{ex7-1}
 \dot x(t) = A_{-1}\dot x(t-1) + A_0x(t),\hskip10mm
A_{-1}=\begin{pmatrix} 1& 1\cr0& 1\end{pmatrix}, \quad  A_0=
 \begin{pmatrix} \alpha& 0\cr 0& \beta\end{pmatrix}.
\end{equation}
One can check that the characteristic equation is $\det \Delta
(\lambda)= (\alpha -\lambda + \lambda e^{-\lambda})(\beta -
\lambda + \lambda e^{-\lambda}) =0$ and for $\alpha\neq\beta$
there are two sequences of eigenvectors, such that
$||v^1_n-v^2_n|| \to 0,$ as $n\to \infty$. It is clear that such  family 
vectors do not form a Riesz basis.}
\end{example}
\section{Stability}
By stability we mean here asymptotic stability. For our neutral type system, as for several infinite dimensional systems, we have essentially two notions of asymptotic stability : exponential (or uniform) stability and strong stability. 
\begin{defn}
A linear system in a Banach space $\mathcal X$ is exponentially stable if the $\e^{\A t}$ semigroup  verifies: 
$
\exists M_{\omega} >1, \quad \exists \omega >0, \quad \forall x, \quad  \Vert \e^{\A t}x \Vert \le M_\omega \e^{-\omega t} \Vert x\Vert.$
 The system is strongly stable if
$
\forall x, \quad \Vert \e^{\A t}x \Vert \to 0, \quad \mathrm{as} \quad t \to \infty.
$
\end{defn}
The problem of exponential stability was widely described in several classical works. An sufficiently exhaustive analysis may be found in \cite{neerven} (see also the references therein and the bibliographic notes). In our case the exponential stability is completely determinated by the spectrum of the operator ${\cal A}$. It is a well known result for some linear neutral type systems: the spectrum has to be bounded away
from the imaginary axis (cf.  \cite[Theorem 6.1]{henry74}).
\begin{thm}\label{thmES0}
The system (\ref{2s}) is
exponentially  stable if and only if    $ \sigma ({\mathcal
A})\subset \{\lambda : \Re \lambda \le -\alpha <0\}$.
\end{thm}
We can partially
reformulate in terms of the matrix $A_{-1}$ the condition on
the spectrum  $\sigma ({\mathcal A})$. 
\begin{thm}\label{thmES}
System (\ref{2s}) is
exponentially  stable if and only if the following conditions are verified \par $i)\quad \sigma ({\mathcal
A})\subset \{\lambda : \Re \lambda < 0\}$ \par
$ii)\quad
\sigma (A_{-1})\subset \{\lambda : |\lambda| < 1\}.$
\end{thm}
It can be interesting how the condition ii) of Theorem \ref{thmES} may be formulated for the case of a general linear operator $K$ in the system \eqref{1}.

We would like to study more deeply the problem of asymptotic non\-exponen\-tial stability. To this end, we recall some important abstract result in this domain. We have the following
\begin{thm}\label{th3}
 Let $e^{At}$, $t\ge0$ be a $C_0$-semigroup in the Banach space $X$ and $A$ be the infinitesimal generator of the semigroup. Assume that  $(\sigma(\A)\cap(i{\Rset}))$ is at most countable and the operator $\A^*$ has no pure imaginary eigenvalues. Then $e^{\A t}$ is strongly asymptotically stable (i.e. $e^{\A t}x\to0$,
$t\to+\infty$ as $x\in X$) if and only if one of the
following conditions is valid:
\begin{itemize}
\item[i)] There exists a norm $\|\cdot\|_1$, equivalent
to the initial one $\|\cdot\|$, such that the semigroup
$e^{At}$  is contractive according to this norm:
$\|e^{At}x\|_1\le\|x\|_1$, $\forall x\in X$, $t\ge0$;

\item[ii)] The semigroup $e^{At}$ is uniformly bounded:
$\exists C>0$ such that $\|e^{At}\|\le C$, $t\ge0$.
\end{itemize}
\end{thm}
The Theorem \ref{th3} was obtained initially in \cite{skl-sh} for a bounded operator $\A$. The main idea were later used in \cite{lyub} for the case of unbounded operator ${\A}$, see also \cite{arendt}  for another  approach. The proof in \cite{lyub} follow the scheme of the first result \cite{skl-sh}. The
development of this theory concerns a large class of
differential equations in Banach space (see \cite{neerven}
and references therein). A more genral result on the asymtotic behavior of the semigroup with respect to an arbitrary asymptotic was recently obtained in \cite{sklyar}.

Our main result on the asymptotic stability of the neutral type system \eqref{1}--\eqref{2s} is the following one.
\begin{thm}\label{thm-stab}     
Assume $\sigma ({\mathcal A})\subset \{\lambda : \Re
\lambda < 0\}$ and  $\sigma_1\equiv \sigma
 (A_{-1})\cap \{\lambda : |\lambda |=1\} \ne \emptyset$ Then the following three mutually
exclusive possibilities exist:
\begin{enumerate}
\item[i)] the part of the spectrum $\sigma_1$ consists of simple
eigenvalues only, i.e. to each eigenvalue corresponds a
one-dimensional 
eigenspace and there are no rootvectors. In this case
system (\ref{2s}) is asymptotically stable.

\item[ii)] the matrix $A_{-1}$ has a Jordan
block, corresponding to  $\mu\in \sigma_1.$ In this case
system (\ref{2s}) is unstable.

\item[iii)] there are no Jordan blocks, corresponding to
eigenvalues in $\sigma_1,$ but there exists  $\mu\in
\sigma_1$ whose eigenspace is at least two-dimensional. In
this case system (\ref{2s}) can be stable as well as
unstable. Moreover, there exist two systems with the same
spectrum, such that one of them is stable while the other one
is unstable.
\end{enumerate}
\end{thm}
The last case may be illustrated by a non trivial example (see also \cite{rsr-stab} for an example given partially in the $M_2$--space framework).
\begin{example}{\em  (Rabah-Sklyar-Barkhaev \cite{rsb-Rep})
Consider the system
$$
\dot z(t) - A_{-1}\dot z(t-1)=A_0z(t)
$$
with
$$
A_{-1}=\begin{pmatrix}-1& 0 \cr 0 & -1\end{pmatrix}, \qquad A_{0}=\begin{pmatrix}-1& \gamma \cr 0 & -1\end{pmatrix}, \qquad \gamma=0 \quad \mbox{or}\quad 1.
$$ 
We have: $\sigma(\A) = \{\lambda : \lambda \e^\lambda + \lambda+ \e^\lambda=0\} $ in $\Cset^-$, this can be proved by Pontriaguin Theorem \cite{pontriag}. The multiplicity of eigenvalues is clearly 2, and they do not depend of $\gamma$. The system is stable for $\gamma=0$ and unstable for $\gamma \ne 0$.
}
\end{example}
\section{Stabilizability}
We say that the system (2)
 is stabilizable if there exists a linear feedback control $u(t)=F(z_t(\cdot))=F(z(t+\cdot))$ such
that the system (2) becomes asymptotically stable.

It is obvious that for linear systems in finite dimensional spaces
the linearity of the feedback implies that the control is bounded
in every neighbourhood of the origin. For infinite dimensional
spaces the situation is much more complicated. The boundedness of
the feedback law $u=F(z_t(\cdot))$ depends on the topology of the
state space.

When the asymptotic
 stabilizability is achieved by a feedback law which does not change the state
 space and is bounded with respect to  the topology of the state space, then we
 call it {\it regular} asymptotic stabilizability. Under our assumption on the state
 space, namely $H^1([-1,0],\Cset^n),$ the natural linear feedback is
 \begin{equation}\label{f1}
 Fz(t+\cdot)=\int_{-1}^0F_2(\theta)\dot z(t+\theta)\dd t +\int_{-1}^0F_3(\theta)
 z(t+\theta)\dd t,
  \end{equation}
  where $F_2(\cdot),F_3(\cdot)\in L_2(-1,0;\Cset^n).$

   Several authors (see for example \cite{hale-lunel,tarn,pand2,dusser} and references therein) use feedback
   laws which for our system may take the form
   \begin{equation}\label{f2} 
   \sum_{i=1}^kF_i \dot z(t-h_i)+\int_{-
   1}^0F_2(\theta)\dot z(t+\theta)\dd t +\int_{-1}^0F_3(\theta) z(t+\theta)\dd t .
   \end{equation} 
   This feedback law is not bounded in $H^1([-1,0],\Cset^n)$ and then
   stabilizability is not regular. As a counterpart, they obtain exponentially stable closed loop systems. If the original system is not formally stable (see \cite{loiseau}), i.e. the pure neutral part (when $A_2=A_3=0$) is not stable, the non regular feedback (\ref{f1}) is necessary to stabilize.
From the operator point of view, the regular feedback law
(\ref{f1}) means a perturbation of the infinitesimal generator
$\A$ by the operator $\B\F$ which is relatively
 $\A$-bounded (cf. \cite{kato}) and verifies $\DD(\A)=\DD(\A+\B\F).$ Such a perturbation does
 not mean, in general, that $\A+\B\F$ is the infinitesimal generator of a $C_0$-semigroup.
 However, in our case, this fact is verified directly \cite{rsr-stab,RSsiam} since
 after the feedback we get also a neutral type system like (\ref{1}) with $\DD(\A)=\DD(\A+\B\F)$ (see below for more details).

 From a physical point of view, $\A$-boundedness of the stabilizing feedback $\F$
 means that the energy added by the feedback remains uniformly bounded in every neighbourhood of 0 (see also another point of view in \cite{loiseau}).
 Hence the problem of {\it regular} asymptotic stabilizability for the systems (\ref{1}),(\ref{2s})
 is to find a linear relatively $\A$-bounded feedback $u=\F x$ such that the operator
 $\A+\B\F$ generates a $C_0$-semigroup ${\mathrm e}^{(\A+\B\F)t}$
 with $\DD(\A+\B\F)=\DD(\A)$ and for which
 $\Vert {\mathrm e}^{(\A+\B\F)t}x \Vert\to 0$, as $t\to\infty$ for all $x\in {\DD(\A)}$. 
The main contribution of this paper is that under some
controllability
  conditions on the unstable poles of the system, we can assign  arbitrarily
  the eigenvalues of the closed loop system into circles centered at the unstable
  eigenvalues of the operator ${\mathcal A}$ with radii  $r_k$ such that $\sum r_k^2<\infty$.
  This is, in some sense, a generalization of the classical pole assignment problem in finite dimensional
  space. 
Precisely we have the following
\begin{thm}\label{thm-stbz}
Consider the system (\ref{1}) under the following assumptions:
\begin{description}
    \item[1)] All the eigenvalues of the matrix $A_{-1}$ satisfy $|\mu
    |\le 1$.
    \item[2)]  All the eigenvalues $\mu_j\in\sigma_1\stackrel{\mathrm{def}}{=}\sigma ( A_{-1})\cap \{
    z : |z|=1 \}$ are simple (we denote their index $j\in I$).
\end{description}

 Then the system (\ref{1})
  is regularly asymptotic stabilizable if 
 \begin{description}
  \item[3)] $\mathrm{rank} \begin{pmatrix}
  \Delta_{\A}(\lambda)&B \end{pmatrix} = n $ for all $\Re\lambda \ge 0$,
  where $$\Delta_{\A}(\lambda)= -\lambda I +
\lambda {\mathrm e}^{-\lambda}A_{-1} + \lambda \int^0_{-1}{\mathrm
e}^{\lambda s} A_2(s){\mathrm d}s  + \int^0_{-1}{\mathrm
e}^{\lambda s} A_3(s){\mathrm d}s,$$
  \item[4)] $\mathrm{rank} \begin{pmatrix}\mu I-A_{-1}& B\end{pmatrix}=n$ for all $\vert\mu\vert=1$.
\end{description}
\end{thm}
\section{Exact Controllability}
The problem of controllability for delay systems was considered by several authors
in different framework.  One approach
is based on the analysis of time delay system in a
module framework (space over ring, see \cite{morse}). In this case the controllability
problem is considered in
a formal way using different interpretations of the Kalman rank condition.
Another approach is based on the analysis of time delay systems in vector spaces with
finite or infinite dimension. A powerful tool is to consider a delay system as
a system in a Banach functional space, this approach was developed
widely in \cite{hale}. Because the state space for delay systems is a functional space,
the most important notion is the function space controllability.
A first important contribution
in the characterization of null functional controllability was given by Olbrot \cite{olbrot} 
by using
 some finite dimensional tools as $(A,B)$-invariant
subspaces for an extended system.
For retarded systems one can refer to \cite{manitius} (and references therein)
 for the analysis of function space controllability in abstract Banach spaces.
 The case of neutral type systems with discrete delay was also considered in such a
 framework (see O'Connor and Tarn \cite{tarn} and references therein).  A general analysis
 of the  time delay systems in infinite dimensional spaces is given in the book
 \cite{bensoussan} where several
 methods and references are given.

The problem considered in this paper is close to that studied in
\cite{tarn}. 
In this  work the exact controllability problem was considered
for neutral type systems with discrete delay using a  semigroup approach in
Sobolev spaces $W_2^{(1)}$ and a boundary control problem.

We
consider the problem of controllability for distributed delay system of neutral type
in the space $M_2(-h,0; \Cset^n)=\Cset^n \times L_2(-h,0; \Cset^n)$
which is natural for control problems.

The semigroup theory developed here is
based on the Hilbert space model introduced in \cite{burns}.
One of our result is a generalization of the result
in \cite{tarn}. The main non trivial precision is the time of controllability.
We generalize the results given \cite{jacob-lang} for the case of a single input
and one localized delay (see also \cite{banks,rodas}).
The approach developed here
 is different from that of
\cite{tarn}. Our main results are based on the characterization of controllability as
a moment problem and using some recent results on  the solvability of this problem
(see \cite{avdo} for the main tools used here). Using a precise Riesz basis in the space
$M_2(-h,0; \Cset^n)$ we can
give a characterization of
null-controllability and of the minimal time of controllability.

The reachability set at time $T$  is defined by
$$\RR_T=\left\{\int_0^Te^{\A t}\B u(t)dt : u(\cdot)\in L_2(0,T;\Cset^n)
\right\}$$

It is easy to show that $\RR_{T_1}\subset \RR_{T_2}$ as $T_1<T_2.$ An important
 result
is that $\RR_T \subset \DD(\A) \subset M_2$. This non-trivial fact permits to formulate
the null-controllability problem in the following setting:\\[.5ex]
{\it  i) To find maximal possible set $\RR_T$ (depending on $T$)};\\[.5ex]
{\it ii) To find minimal $T$ for which the set $\RR_T$ becomes maximal
possible, i.e. $\RR_T = \DD(\A)$}.
\begin{defn}
The system (\ref{2s}) is said null-controllable at the time $T$ if $\RR_T = \DD(\A)$
\end{defn}
The main tool is to consider the null-controllability problem as a problem of moments.
\subsection{The moment problem}
In order to formulate the moment problem we need a Riesz basis in the Hilbert space $M_2$.
We recall that a Riesz basis is a basis which may be transformed to an orthogonal
basis with respect to another equivalent scalar product. Each Riesz basis possesses
a biorthogonal basis. Let $\{\varphi\}$ be a Riesz basis in $M_2$ and $\{\psi\}$ 
the corresponding
biorthogonal basis. Then for each $x \in M_2$ we have
$
x=\sum_{\varphi \in\{\varphi\}}\langle x, \psi \rangle \varphi.
$
In a separable Hilbert space there always exists a Riesz basis.

A state $x=\begin{pmatrix}y \cr z(\cdot)\end{pmatrix} \in M_2$ is reachable at time $T$ by a control
 $u(\cdot) \in L_2(0,T;\Cset^r)$ iff the steering condition
\begin{equation}\label{1.1}
x=\begin{pmatrix}y\cr z(\cdot)\end{pmatrix}=\int_0^T\e^{{\mathcal A}t}{\mathcal B}u(t)\dd t.
\end{equation}
holds. This steering condition may be expanded using the basis $\{\varphi\}$. A state $x$ is
reachable iff
$$
\sum_{\varphi \in\{\varphi\}}\langle x, \psi\rangle \varphi= \sum_{\varphi
\in\{\varphi\}}\int_0^T\langle \e^{{\mathcal A}t}{\mathcal B}u(t), \psi\rangle\dd t \varphi,
$$
for some $u(\cdot) \in L_2(-h,0; \Rset^r)$.
Then  the steering condition (\ref{1.1}) can be substituted by the following system of
equalities
\begin{equation}\label{2.1}
\langle x, \psi\rangle = \int_0^T\langle \e^{{\mathcal A}t}{\mathcal B}u(t), \psi\rangle\dd t ,
\quad \psi \in\{\psi\}.
\end{equation}
Let $\{b_1, \dots, b_r \}$ be an arbitrary basis in ${\mathrm {Im}}B$, the image
of the matrix $B$ and ${\mathbf b}_i=\begin{pmatrix}b_i\cr 0\end{pmatrix} \in M_2, \ i=1, \dots,r$.
 Then the right hand side of (\ref{2.1}) takes the form
\begin{equation}
  \label{2.2}
\int_0^T\langle \e^{{\mathcal A}t}{\mathcal B}u(t), \psi\rangle\dd t
=\sum_{i=1}^r\int_0^T\langle \e^{{\mathcal A}t}\bfb_i,\psi\rangle u_i(t)\dd t.
\end{equation}
Effectiveness of the proposed approach becomes obvious if we assume that the
operator $\mathcal{A}$ possess a Riesz basis of eigenvector. This situation is
characteristic, for example, for control systems of hyperbolic type when
$\mathcal{A}$ is skew-adjoint ($\mathcal{A}^* = -\mathcal{A}$) and has a
compact resolvent (see, for example, [1], [16], [17]). Let in this case
$\{\varphi_k\},$ $k\in\mathbb{N},$ be a orthonormal eigenbasis with
$\mathcal{A}\varphi_k = i\lambda_k\varphi_k,$ $\lambda_k\in\mathbb{R}.$
Assuming for simplicity $r = 1,$ $b_1 = b = \sum_k\alpha_k\varphi_k,$ $\alpha_k
\neq 0,$ we have from (4), (5)
\begin{equation}
  \label{2.3}
\frac{x_k}{\alpha_k}=\int_0^Te^{-i\lambda_kt}u(t)dt, \quad k\in\mathbb{N},
\end{equation}
where $x = \sum_k x_k\varphi_k.$ Equalities (6) are a non-Fourier trigonometric
moment problem whose solvability is closely connected with the property for the
family of exponentials $e^{-i\lambda_kt},$ $k\in\mathbb{N},$ to form a Riesz
basis on the interval $[0,T]$ ([1]). In particular, if $e^{-i\lambda_kt}$ forms
a Riesz basis of $L_2[0,T_0]$ then one has
\begin{equation}
  \label{2.4}
\mathcal{R}_T =
\left\{x:\sum_k\left(\frac{x_k}{\alpha_k}\right)^2<\infty\right\} \quad
\hbox{for all } T \geq T_0.
\end{equation}
Obviously formula (\ref{2.4}) gives the complete answer to the both items of the
controllability problem. Returning now to neutral type systems we observe that
the operator $\mathcal{A}$ given in is not skew-adjoint and, moreover, does
not possess a basis even of generalized eigenvectors. So the choice of a proper
Riesz basis to transform the steering condition in a moment problem is an essentially more complicated
problem.
\subsection{The choice of basis}
In order to design the needed basis for our case we use  spectral the 
properties of the operator $\A$ obtained in \cite{rsr-stab}.
Let $\mu_1,\ldots,\mu_\ell$, \ $\mu_i\not=\mu_j$ be eigenvalues of
$A_{-1}$ and let the integers $p_m$ be defined as :
$
{\rm dim}\,(A_{-1}-\mu_mI)^n=p_m,\ m=1,\ldots,\ell.
$
Denote by
$$
\lambda_m^{(k)}=\frac{1}{h}\left (\ln|\mu_m|+{\mathrm i}(\arg\mu_m+2\pi k)\right ),\
m=1,\ldots,\ell;\  k\in\Zset,
$$
and let $L_m^{(k)}$ be the circles of the fixed radius $r\le
r_0=\frac13\min|\lambda_m^{(k)}- \lambda_i^{(j)}|$ centered at
$\lambda_m^{(k)}$.

Let $\{V_m^{(k)}\}_{\frac{k\in\,\Zset\hfill}{m=1,\ldots,\ell}}$ be a
family of $\A$-invariant subspaces given by
$$V_m^{(k)}=P_m^{(k)}M_2,\qquad P_m^{(k)}=\frac{1}{2\pi {\mathrm i}}\int_{L_m^{(k)}}R(\A,\lambda)d\lambda.$$
The following theorem plays an  essential role in our approach
\begin{thm}\cite{rsr-CR,rsr-stab}
There exists $N_0$ large
enough such that for any $N\ge N_0$\\
i) ${\rm dim}\,V_m^{(k)}=p_m$, $k\ge N$,\\
ii) the family $\{V_m^{(k)}\}_{\frac{|k|\ge
N\hfill}{m=1,\ldots,\ell}}\cup \widehat V_N$ forms a Riesz
basis (of subspaces) in $M_2$, where $\widehat V_N$ is a
finite-dimensional subspace (${\rm dim}\,\widehat
V_N=2(N+1)n$) spanned by all generalized eigenvectors
corresponding to all eigenvalues of $\A$ located outside of all
circles $L_m^{(k)}$, $ |k|\ge N$, $m=1,\ldots,\ell$.
\end{thm}
Using this theorem we construct a Riesz basis $\{\varphi\}$ of the form
$$
\left \{\varphi_{m,j}^k,
|k|>N;
m=1,\ldots,l;
j=1,\ldots,p_m\right \}
\cup\left \{\hat{\varphi}_j^N, j=1,\ldots,2(N+1)n\right \}
$$
where for any $m=1,\ldots,l,$ and $k:|k|>N$ the collection
$\{\varphi_{m,j}^k\}_{j=1,\ldots,p_m}$ is in a special way chosen basis of
$V_m^{(k)}$ and $\{\hat{\varphi}_j^N\}_{j=1,\ldots,2(N+1)n}$ is a basis of
$\hat{V}_N.$ In this basis equalities (4) with regard to (5) turns into a
moment problem with respect to a special collection of quasipolynomials.
Analyzing the mentioned moment problem by means of the methods given in [1] we
obtain our main results concerning the null-controllability problem.
\section{The main results}
The characterization  of the null-controllability is given by the following Theorem.
\begin{thm}\label{thm4}
The system (\ref{2s}) is null-controllable by controls from
$L_2(0,T)$ for some $T>0$ iff the following two conditions hold:\\[0.ex]
i)
$
\rank[\begin{matrix}\Delta_{\mathcal A}(\lambda)& B\end{matrix}]=n, \quad \forall \lambda \in \Cset;
$
where
\begin{eqnarray*}
\Delta_{\mathcal A}(\lambda)= -\lambda I
+ \lambda \e^{-\lambda h }A_{-1}
 + \lambda
\int^0_{-h}\e^{\lambda s} A_2(s){\mathrm d}s
 + \int^0_{-h}\e^{\lambda
s} A_3(s){\mathrm d}s.
\end{eqnarray*}
ii) $\rank [\begin{matrix}B&A_{-1}B&\cdots&A_{-1}^{n-1}B\end{matrix}]=n$.

\end{thm}
The main results on the time of controllability are as follows.
 \begin{thm} \label{thm5}
Let the conditions i) and ii) of Theorem~\ref{thm4}
 hold. Then
 \begin{enumerate}
\item[i)] The system (\ref{2s}) is null-controllable at the time $T$ as $T>nh$;
\item[ii)] If the system (\ref{2s}) is of single control ($r=1$),
then the estimation of the time of controllability in i) is exact,
i.e. the system is not controllable at time $T=nh$.
\end{enumerate}
 \end{thm}
 For the multivariable case, the time depends on some controllability indices. suppose that
 $\dim B =r$. Let $\{b_1, \dots, b_r\}$ be an arbitrary basis noted $\beta$.
 Let us introduce a set integers.
 We denote by $B_i=\begin{pmatrix}{b_{i+1}, \dots,b_r}\end{pmatrix}, \ i=0, 1, \dots, r-1$, which
 gives in particular $B_0=B$ and $B_{r-1}={b_r}$ and we put
 formally $B_r=0$.
Let us consider the integers 
$$
n_i^\beta=\rank[\begin{matrix}B_{i-1}&A_{-1}B_{i-1}&\cdots&A_{-1}^{n-1}B_{i-1}\end{matrix}], \quad i=1, \dots,r,
$$
corresponding to the basis $\beta$. We need in fact the integers
$$
m_i^\beta=n_{i-1}^\beta-n_{i}^\beta, 
$$
Let us denote by
$$
m_{\mathrm {min}}=\max\limits_{\beta}m_1^\beta \qquad
m_{\mathrm {max}} = \min\limits_{\beta}\max\limits_{i}m_i^\beta ,
$$
for all possible choice of a basis  $\beta$.

 The main result for the multivariable case is the following Theorem.
 \begin{thm}\label{thm6}
 Let the conditions i) and ii) of the Theorem~\ref{thm4} hold, then
 \begin{enumerate}
\item[i)] The system (\ref{2s}) is null-controllable at the time $T>m_{\mathrm {max}}h$;
\item[ii)] The system (\ref{2s}) is not null-controllable at the time $T< m_{\mathrm {min}}h$.
\end{enumerate}
\end{thm}

 The proofs are based on the construction of a special Riesz basis of
 $\A$-invariant subspaces in the space $M_2$
 according to
 \cite{rsr-CR,rsr-stab} and on the analysis of the properties of some quasi-exponential
 functions to be a Riesz basis in $L_2(0,T)$ depending of the time $T$ \cite{avdo}.

\end{document}